\documentclass{amsart}
\usepackage[dvips]{graphicx}
\usepackage{amscd}
\usepackage{amsmath}
\usepackage{amsxtra}
\usepackage{amsfonts}
\usepackage{amssymb}
\newtheorem{theorem}{Theorem}[section]

\newtheorem{proposition}[theorem]{Proposition}
\theoremstyle{definition}
\newtheorem{definition}[theorem]{Definition}
\newtheorem{remark}[theorem]{Remark}

\newtheorem{example}[theorem]{Example}
\theoremstyle{remark}

\renewcommand{\theclaim}{\textup{\theclaim}}

\numberwithin{equation}{section}

\def\openone

{\mathchoice

{\hbox{\upshape \small1\kern-3.3pt\normalsize1}}

{\hbox{\upshape \small1\kern-3.3pt\normalsize1}}

{\hbox{\upshape \tiny1\kern-2.3pt\SMALL1}}

{\hbox{\upshape \Tiny1\kern-2pt\tiny1}}}

\makeatletter

\newbox\ipbox

\newcommand{\diracb}[1]{\left\langle #1\mathrel{\mathchoice

{\setbox\ipbox=\hbox{$\displaystyle \left\langle\mathstrut #1\right.$}

\vrule height\ht\ipbox width0.25pt depth\dp\ipbox}

{\setbox\ipbox=\hbox{$\textstyle \left\langle\mathstrut #1\right.$}

\vrule height\ht\ipbox width0.25pt depth\dp\ipbox}

{\setbox\ipbox=\hbox{$\scriptstyle \left\langle\mathstrut #1\right.$}

\vrule height\ht\ipbox width0.25pt depth\dp\ipbox}

{\setbox\ipbox=\hbox{$\scriptscriptstyle \left\langle\mathstrut #1\right.$}

\vrule height\ht\ipbox width0.25pt depth\dp\ipbox}

}\right. }

\newcommand{\dirack}[1]{\left. \mathrel{\mathchoice

{\setbox\ipbox=\hbox{$\displaystyle \left.\mathstrut #1\right\rangle$}

\vrule height\ht\ipbox width0.25pt depth\dp\ipbox}

{\setbox\ipbox=\hbox{$\textstyle \left.\mathstrut #1\right\rangle$}

\vrule height\ht\ipbox width0.25pt depth\dp\ipbox}

{\setbox\ipbox=\hbox{$\scriptstyle \left.\mathstrut #1\right\rangle$}

\vrule height\ht\ipbox width0.25pt depth\dp\ipbox}

{\setbox\ipbox=\hbox{$\scriptscriptstyle \left.\mathstrut #1\right\rangle$}

\vrule height\ht\ipbox width0.25pt depth\dp\ipbox}

} #1\right\rangle}

\newcommand{\ltwor}{L^{2}\left(\mathbb{R}\right)}

\newcommand{\ltworn}{L^{2}\left(\mathbb{R}^n\right)}

\newcommand{\ltwozn}{l^{2}\left(\mathbb{Z}^n\right)}

\newcommand{\Trace}{\operatorname*{Trace}}

\newcommand{\Tper}{\mathcal{T}_{per}}

\newcommand{\Rn}{\mathbb{R}^n}

\newcommand{\Zn}{\mathbb{Z}^n}

\input cyracc.def

\begin{document}
\title[Multiwavelets and multiscaling functions]{Some equations relating multiwavelets and multiscaling functions}
\author{Dorin Ervin Dutkay}
\address{Department of Mathematics\\
The University of Iowa\\
14 MacLean Hall\\
Iowa City, IA 52242-1419\\
U.S.A.}
\email{ddutkay@math.uiowa.edu}
\thanks{}

\begin{abstract}
The local trace function introduced in \cite{Dut} is used to derive equations that relate multiwavelets and multiscaling functions in  the context of
 a generalized multiresolution analysis, without appealing to filters. A construction of normalized tight frame wavelets is given.
Particular instances of the construction include normalized tight frame and orthonormal wavelet sets.  
\end{abstract} 

\maketitle 
\tableofcontents

\section{\label{introduction}Introduction}
A wavelet is a function $\psi\in\ltwor$ such that 
$$\{D^jT_k\psi\,|\,j\in\mathbb{Z},k\in\mathbb{Z}\}$$
is an orthonormal basis for $\ltwor$, where 
$$Df(\xi)=\sqrt{2}f(2\xi),\quad T_kf(\xi)=f(\xi-k),\quad(\xi\in\mathbb{R},f\in\ltwor,k\in\mathbb{Z}).$$
Many examples of wavelets have been produced using the concept of multiresolution analysis (MRA) (see \cite{Dau}). A MRA 
is a nest of subspaces $(V_n)_{n\in\mathbb{Z}}$ of $\ltwor$ with the following properties:
\begin{equation}\label{eq0_1}
V_n\subset V_{n+1},\quad(n\in\mathbb{Z});
\end{equation}
\begin{equation}\label{eq0_2}
f\in V_n\mbox{ iff }Df\in V_{n+1};
\end{equation}
\begin{equation}\label{eq0_3}
\overline{\bigcup_{n\in\mathbb{Z}}V_n}=\ltwor;
\end{equation}
\begin{equation}\label{eq0_4}
\bigcap_{n\in\mathbb{Z}}V_n=\{0\};
\end{equation}
\begin{equation}\label{eq0_5}
\mbox{ There exists }\varphi\in V_0\mbox{ such that }\{T_k\varphi\,|\, k\in\mathbb{Z}\}\mbox{ is an orthonormal basis for }V_0.
\end{equation}
$\varphi$ is called a scaling function.
\par
To construct wavelets, one has to find functions $\psi$ such that 
$\{T_k\psi\,|\,k\in\mathbb{Z}\}$ is an orthonormal basis for $W_0:=V_1\ominus V_0$. 
\par
Many examples, due to Journe and others (\cite{DL}), show that there are wavelets which are not associated to MRAs. The theory 
developed by Baggett (\cite{BMM},\cite{BM}) shows that every orthogonal wavelet is associated to a similar, more general structure 
called generalized multiresolution analysis (GMRA) which satisfies the conditions (\ref{eq0_1})-(\ref{eq0_4}) while condition 
(\ref{eq0_5}) is replaced by a weaker one:
\begin{equation}\label{eq0_6}
V_0\mbox{ is invariant under all integer trasnlations }T_k.
\end{equation}
In the context of a MRA, wavelets are constructed from scaling functions via filters (\cite{Dau}). In the GMRA case, one can
still get some scaling functions in $V_0$, namely, there are functions $\phi_1,...,\phi_n,...\in V_0$ such that 
$$\{T_k\phi_i\,|\,i\in\mathbb{N},k\in\mathbb{Z}\}$$
is a normalized tight frame for $V_0$. 
\par
We recall that a set of vectors $\{e_i\,|\, i\in I\}$ in a Hilbert space $H$ is a frame if there are some positive constants
$A,B>0$ such that
$$A\|f\|^2\leq\sum_{i\in I}\left|\left\langle f\,|\, e_i \right\rangle\right|^2\leq B\|f\|^2,\quad (f\in H).$$
If $A=B=1$ it is called a normalized tight frame (NTF). 
\par
In the GMRA situation, the wavelets can be again constructed using filters but substantial complications appear because,
instead of just one filter, as it was in the case of a MRA, now one has to use a matrix of filters. 
\par
In this paper we analyse the relation between scaling functions and wavelets without the use of filters. This relation is described 
in three theorems:
\par
1. In theorem \ref{th1_1} we assume that the scaling functions are given and offer necessary and sufficient conditions for a set of functions 
to be an associated wavelet. 
\par
2. In theorem \ref{th1_2} we start with a wavelet and derive equations that characterize the associated scaling functions. 
\par
3. In theorem \ref{th1_5} we show that if two sets of functions are related by equations similar to those that link scaling functions 
and wavelets then one of the sets will be indeed a wavelet. (However the other set is not necessarily the corresponding scaling function.)
\par
In section \ref{some} we list some definitions and results in preparation for the main part, which is section \ref{main} where 
the results are proved. 
In section \ref{construction}, we describe a general procedure for constructing normalized tight frame wavelets. All wavelet sets and 
normalized tight frame wavelet sets can be obtained with this procedure provided the initial data is chosen appropriately. We 
end with an example of a NTF wavelet which has a piecewise linear square in the Fourier domain. 

\section{\label{some} Some definitions and tools} 
Throughout the paper we will work with an $n\times n$ dilation matrix $A$ which preserves the lattice $\Zn$, that is: 
all the eigenvalues $\lambda$ of $A$ have $|\lambda|>1$ and $A\Zn\subset\Zn$. Define the translation and dilation operators on 
$\ltworn$:
$$T_kf(\xi)=f(\xi-k),\quad D_Af(\xi)=|\operatorname*{det}A|^{\frac{1}{2}}f(A\xi),\quad(\xi\in\Rn,f\in\ltworn,k\in\Zn).$$
For a subset $\Psi$ of $\ltworn$ define the affine system 
$$X(\Psi):=\{D_A^jT_k\psi\,|\,j\in\mathbb{Z},k\in\Zn,\psi\in\Psi\}.$$
$\Psi$ is called a normalized tight frame (orthogonal) multiwavelet if $X(\Psi)$ is a normalized tight frame (orthonormal basis) for 
$\ltworn$. 
\par
A generalized multiresolution analysis (GMRA) is a nest of closed subspaces $(V_n)_{n\in\mathbb{Z}}$ of $\ltworn$ with the following properties:
\begin{equation}\label{eq00_1}
V_n\subset V_{n+1},\quad(n\in\mathbb{Z});
\end{equation}
\begin{equation}\label{eq00_2}
f\in V_n\mbox{ iff }D_Af\in V_{n+1};
\end{equation}
\begin{equation}\label{eq00_3}
\overline{\bigcup_{n\in\mathbb{Z}}V_n}=\ltworn;
\end{equation}
\begin{equation}\label{eq00_4}
\bigcap_{n\in\mathbb{Z}}V_n=\{0\};
\end{equation}
\begin{equation}\label{eq00_5}
T_kV_0=V_0,\quad(k\in\Zn).
\end{equation}
\par
A multiscaling function associated to a GMRA is a subset $\Phi$ of $\ltworn$ such that 
$\{T_k\varphi\,|\,k\in\Zn,\varphi\in\Phi\}$ is a normalized tight frame for $V_0$.
\par
The Fourier transform is given by
$$\widehat{f}(\xi)=\int_{\Rn }f(x)e^{-i\left\langle x\,|\,\xi\right\rangle}\,dx,\quad(\xi\in\Rn ).$$
If $V$ is closed subspace of a Hilbert space $H$ and $f\in H$ we denote by $P_V$ the projection onto $V$ and
by $P_f$ the operator defined by:
$$P_f(v)=\left\langle v\,|\,f\right\rangle f,\quad(v\in H).$$  
\par
The main tool needed for our analysis will be the local trace function introduced in \cite{Dut}. For details, several properties and
the appropriate references 
we refer the reader to that paper. We recall below the definition and some properties that will be used here. 
The local trace function is associated to shift invariant spaces. 
\begin{definition}\label{def00_1}
A closed subspace $V$ of $\ltwor$ is called shift invariant (or shortly SI) if 
$$T_kV=V,\quad(k\in\Zn).$$
If $\mathcal{A}$ is a subset of $\ltworn$ then we denote by $S(\mathcal{A})$ the shift invariant space generated
by $\mathcal{A}$,
$$S(\mathcal{A})=\overline{\operatorname*{span}}\{T_k\varphi\,|\,k\in\Zn,\varphi\in\mathcal{A}\}.$$
\end{definition}
\par
\begin{definition}\label{def00_2}
Let $V$ be a shift invariant subspace of $\ltworn$. A subset $\Phi$ of $V$ is called a normalized tight frame
generator (or NTF generator) for $V$ if
$$\{T_k\varphi\,|\,k\in\Zn,\varphi\in\Phi\}$$
is a NTF for $V$.
\end{definition}
\par
Shift invariant spaces have been studied in connection not only to wavelets but also to spline systems, Gabor systems or approximation theory. 
The local trace function is constructed using some fiberization techniques introduced in \cite{H} and developed by A.Ron, Z.Shen,
M. Bownik, Z. Rzeszotnik and others (\cite{RS1}, \cite{RS2}, \cite{RS3}, \cite{Bo1}, \cite{BoRz}). 
These "fiberization" tools include the range function. For more information on the range function we refer to \cite{H},\cite{Bo1} and \cite{Dut}.
The periodic range function is a measurable map from $\Rn$ to the projections (or subspaces) of $\ltwozn$ satisfying the periodicity:
$$J_{per}(\xi+2k\pi)=\lambda(k)^*\left(J_{per}(\xi)\right),\quad(k\in\Zn,\xi\in\Rn),$$
where $\lambda$ denotes the shift on $\ltwozn$,
$$(\lambda(k)\alpha)(l)=\alpha(l-k),\quad(l\in\Zn,k\in\Zn).$$
$\Tper$ is defined on $\ltworn$ by
$$\Tper f(\xi)=(\widehat{f}(\xi+2k\pi))_{k\in\Zn},\quad(\xi\in\Rn,f\in\ltworn).$$
Periodic range functions are associated to shift invariant subspaces in a unique way, the connection being described by the 
following theorem due to Helson:
\begin{theorem}\label{th00_3}
A closed subspace $V$ of $\ltworn$ is shift invariant if and only if 
$$V=\{f\in\ltworn\,|\, \Tper f(\xi)\in J_{per}(\xi)\mbox{ for a.e. }\xi\in\mathbb{R}^n\},$$
for some measurable periodic range function $J_{per}$. The correspondence between $V$ and $J_{per}$ is 
bijective under the convention that range functions are identified if they are equal a.e. Furthermore, if 
$V=S(\mathcal{A})$ for some countable $\mathcal{A}\subset\ltworn$, then 
$$J_{per}(\xi)=\overline{\operatorname*{span}}\{\Tper\varphi(\xi)\,|\,\varphi\in\mathcal{A}\},\quad \mbox{for
a.e. }\xi\in\Rn.$$
\end{theorem}
The local trace function is defined as follows:
\begin{definition}\label{def00_4}
Let $V$ be a SI subspace of $\ltworn$, $T$ a positive operator on $\ltwozn$ and let $J_{per}$ be the range function 
associated to $V$. We define the local trace function associated to $V$ and $T$ as the map from $\mathbb{R}^n$ to
$[0,\infty]$ given by the formula 
$$\tau_{V,T}(\xi)=\Trace\left(TJ_{per}(\xi)\right),\quad(\xi\in\mathbb{R}).$$
We define the restricted local trace function associated to $V$ and a vector $f$ in $\ltwozn$ by 
$$\tau_{V,f}(\xi)=\Trace\left(P_fJ_{per}(\xi)\right)(=\tau_{V,P_f}(\xi)),\quad(\xi\in\mathbb{R}^n).$$
\end{definition} 
\par 
Theorems \ref{th00_5} gives a formula for the computation of the local trace function.
\begin{theorem}\label{th00_5}\cite{Dut}
Let $V$ be a SI subspace of $\ltworn$ and $\Phi\subset V$ a NTF generator for $V$. Then for every positive operator $T$ 
on $\ltwozn$ and any $f\in\ltwozn$, 
\begin{equation}\label{eq00_5_1}
\tau_{V,T}(\xi)=\sum_{\varphi\in\Phi}\left\langle T\Tper\varphi(\xi)\,|\,\Tper\varphi(\xi)\right\rangle,\quad\mbox{for a.e. }\xi\in\mathbb{R}^n;
\end{equation}
\begin{equation}\label{eq00_5_2}
\tau_{V,f}(\xi)=\sum_{\varphi\in\Phi}|\left\langle f\,|\,\Tper\varphi(\xi)\right\rangle|^2,\quad \mbox{ for a.e. }\xi\in\mathbb{R}^n.
\end{equation}
\end{theorem}
We should point out that the equations (\ref{eq00_5_1}) and (\ref{eq00_5_2}) show that the local trace function can be calculated
with {\it any} NTF generator. This is the fact that we will use frequently: the local trace function can be computed in two (ore more) different ways 
and the resulting quantities must be equal. 
\par
The next theorem characterizes the NTF generators for a SI space. 
\begin{theorem}\label{th00_5_1}\cite{Dut}
Let $V$ be a SI subspace of $\ltworn$, $J_{per}$ its periodic range function and $\Phi$ a countable subset of $\ltworn$.
Then following affirmations are
equivalent:
\begin{enumerate}
\item
$\Phi\subset V$ and $\Phi$ is a NTF generator for $V$;
\item
For every $f\in\ltwozn$ 
\begin{equation}\label{eq00_5_1_1}
\sum_{\varphi\in\Phi}|\left\langle f\,|\,\Tper\varphi(\xi)\right\rangle |^2=\|J_{per}(\xi)(f)\|^2(=\tau_{V,f}(\xi)),\quad\mbox{for a.e. }\xi\in\Rn
\end{equation}
\item
For every $0\neq l\in\Zn$ and $\alpha\in\{0,1,i\}$,
\begin{equation}\label{eq00_5_1_2}
\sum_{\varphi\in\Phi}|\widehat{\varphi}(\xi)+\overline{\alpha}\widehat{\varphi}(\xi+2l\pi)|^2=\|J_{per}(\xi)(\delta_0+\alpha\delta_l)\|^2(=\tau_{V,\delta_0+\alpha\delta_l}(\xi)),\quad\mbox{for
a.e. }\xi\in\Rn.
\end{equation}
\end{enumerate}
\end{theorem}

\par
The local trace function contains the dimension function $\mbox{dim}_V$ and the spectral function
$\sigma_V$ introduced in \cite{BoRz}. More precisely:
$$\mbox{dim}_V=\tau_{V,I},\quad\sigma_V=\tau_{V,\delta_0}.$$
where 
$$\delta_k(l)=\left\{\begin{array}{ccc}
1&\mbox{if}&k=l,\\
0& &\mbox{otherwise.}
\end{array}\right.
$$
Here are some properties of the local trace function:
\begin{proposition}\label{prop00_6}\cite{Dut}
\begin{enumerate}
\item
If $V_1,V_2$ are orthogonal shift invariant subspaces then
$$\tau_{V_1\oplus V_2,f}=\tau_{V_1,f}+\tau_{V_2,f},\quad(f\in\ltwozn).$$
\item
If $V_1\subset V_2$ are SI subspaces then
$$\tau_{V_1,f}\leq\tau_{V_2,f},\quad(f\in\ltwozn).$$
\item
If $(V_i)_{i\in\mathbb{N}}$ is an increasing set of SI subspaces and $V=\overline{\cup V_i}$, then
$$\tau_{V,f}=\lim_{i\rightarrow\infty}\tau_{V_i,f}.$$
\end{enumerate}
\end{proposition}
\par
The local trace function is well behaved with respect to dilations: the local trace function of the dilation of a SI space can be computed
in terms of the local trace function of the in initial space:
\begin{proposition}\label{prop00_7}\cite{Dut}
Let $V$ be a SI subspace and $A$ an $n\times n$ integer matrix with $\operatorname*{det}A\neq 0$. Then $D_AV$ is
shift invariant and, for every vector $f\in\ltwozn$, 
\begin{equation}\label{eq00_7_1}
\tau_{D_AV,f}(\xi)=\sum_{d\in\mathcal{D}}\tau_{V,D_d^*f}\left(\left(A^*\right)^{-1}(\xi+2d\pi)\right),\quad\mbox{
for
a.e. }\xi\in\mathbb{R}^n,
\end{equation}
where $\mathcal{D}$ is a complete set of representatives of the cosets $\mathbb{Z}^n/A^*\mathbb{Z}^n$
and 
$D_d$ is the linear operator on $\ltwozn$ defined by
$$(D_d\alpha)(k)=\left\{\begin{array}{ccc}
\alpha(l),&\mbox{if}&k=d+A^*l\\
0,& & otherwise
\end{array} 
\right.,\quad(k\in\mathbb{Z}^n,\alpha\in\ltwozn).$$
\end{proposition}
\par
We will also need the following characterization of NTF multiwavelets (see \cite{HW} or \cite{Bo2}).
\begin{theorem}\label{th00_8}
Let $\Psi$ be a finite subset of $\ltworn$. Then $\Psi$ is a NTF multiwavelet iff the following equations are satisfied for a.e. 
$\xi\in\Rn$:
\begin{equation}\label{eq00_8_1}
\sum_{\psi\in\Psi}\sum_{j\in\mathbb{Z}}|\widehat{\psi}|^2((A^*)^j(\xi))=1;
\end{equation}
\begin{equation}\label{eq00_8_2}
\sum_{\psi\in\Psi}\sum_{j\geq0}\widehat{\psi}((A^*)^j(\xi))\overline{\widehat{\psi}}((A^*)^j(\xi+2s\pi))=0,\quad(s\in\Zn\setminus A^*\Zn).
\end{equation}
\end{theorem}
And the last tool that we will need is a relation between multiscaling functions and multiwavelets which was proved in \cite{Dut}.
\begin{theorem}\label{th00_9}
Let $(V_n)_{n\in\mathbb{Z}}$ be a GMRA and $\Psi$ a NTF generator for $W_0:=V_1\ominus V_0$.
Let $\Phi$ be a countable subset
of $\ltworn$. The following affirmations are equivalent:
\begin{enumerate}
\item
$\Phi$ is contained in $V_0$ and is a NTF generator for $V_0$;
\item
The following equations hold: for every $s\in\Zn$,
\begin{equation}\label{eq00_9_1}
\sum_{\psi\in\Psi}\sum_{j\geq 1}\widehat{\psi}((A^*)^j\xi)\overline{\widehat{\psi}}((A^*)^j(\xi+2s\pi))=
\sum_{\varphi\in\Phi}\widehat{\varphi}(\xi)\overline{\widehat{\varphi}}(\xi+2s\pi).
\end{equation}
for a.e. $\xi\in\Rn$.
\end{enumerate}
\end{theorem}

\section{\label{main}Main results}
The first theorem of the section starts with a GMRA for which the scaling functions are given. The theorem characterizes the
wavelets associated to this GMRA. 
\begin{theorem}\label{th1_1}
Let $(V_n)_{n\in\mathbb{Z}}$ be a GMRA and $\Phi$ a NTF generator for $V_0$. Then $\Psi$ is a NTF generator for $W_0:=V_1\ominus V_0$ iff for a.e. $\xi\in\Rn$:
\begin{equation}\label{eq1_1_1}
-\sum_{\varphi\in\Phi}\widehat{\varphi}(\xi)\overline{\widehat{\varphi}}(\xi+2s\pi)=\sum_{\psi\in\Psi}\widehat{\psi}(\xi)\overline{\widehat{\psi}}(\xi+2s\pi),
\quad(s\in\Zn\setminus A^*\Zn);
\end{equation}
$$\sum_{\varphi\in\Phi}\widehat{\varphi}((A^*)^{-1}\xi)\overline{\widehat{\varphi}}((A^*)^{-1}(\xi+2s\pi))-
\sum_{\varphi\in\Phi}\widehat{\varphi}(\xi)\overline{\widehat{\varphi}}(\xi+2s\pi)=$$
\begin{equation}\label{eq1_1_2}
\sum_{\psi\in\Psi}\widehat{\psi}(\xi)\overline{\widehat{\psi}}(\xi+2s\pi),\quad(s\in A^*\Zn).
\end{equation}
\end{theorem}

\begin{proof}
We will use theorem \ref{th00_5_1}. For this, we have to compute $\tau_{W_0,\delta_0+\lambda\delta_s}$ for $s\neq0$, $\lambda\in\{0,1,i\}$. Using the additivity, this will reduce to the computation
of the local trace function for $V_1$. 
\par
Since $V_1=D_AV_0$, the dilation property given in proposition \ref{prop00_7}, yields:
$$\tau_{V_1,\delta_0+\lambda\delta_s}(\xi)=\sum_{d\in\mathcal{D}}\tau_{V_0,D_d^*(\delta_0+\lambda\delta_s)}((A^*)^{-1}(\xi+2d\pi)).$$
\par
The first case we consider is when $s$ is not in $A^*\Zn$. In this case we can assume $0,s\in\mathcal{D}$ and
$$D_d^*(\delta_0+\lambda\delta_s)=\left\{
\begin{array}{ccc}
0&\mbox{if}&d\neq0\mbox{ and }d\neq s\\
\delta_0&\mbox{if}&d=0\\
\lambda\delta_0&\mbox{if}&d=s.
\end{array}\right.$$
So that 
\begin{align*}
\tau_{V_1,\delta_0+\lambda\delta_s}(\xi)&=\tau_{V_0,\delta_0}((A^*)^{-1}\xi)+\tau_{V_0,\lambda\delta_0}((A^*)^{-1}(\xi+2s\pi))\\
&=\sum_{\varphi\in\Phi}|\widehat{\varphi}|^2((A^*)^{-1}\xi)+\sum_{\varphi\in\Phi}|\lambda|^2|\widehat{\varphi}|^2((A^*)^{-1}(\xi+2s\pi)).
\end{align*}
For the last equality we used theorem \ref{th00_5} for $V_0$. 
\par
Also, we have
$$\tau_{V_0,\delta_0+\lambda\delta_s}(\xi)=\sum_{\varphi\in\Phi}|\widehat{\varphi}(\xi)+\overline{\lambda}\overline{\widehat{\varphi}}(\xi+2s\pi)|^2,$$
$$\tau_{W_0,\delta_0+\lambda\delta_s}(\xi)=\tau_{V_1,\delta_0+\lambda\delta_s}(\xi)-\tau_{V_0,\delta_0+\lambda\delta_s}(\xi).$$
and $\Psi$ is a NTF for $W_0$ if and only if for all $s\neq0$

$$\tau_{W_0,\delta_0+\lambda\delta_s}(\xi)=
\sum_{\psi\in\Psi}|\widehat{\psi}(\xi)+\overline{\lambda}\overline{\widehat{\psi}}(\xi+2s\pi)|^2,$$
Hence, if $\Psi$ is a NTF for $W_0$ then: 
\\
If we take $\lambda=0$ then 
\begin{equation}\label{eq1_1_3}
\sum_{\varphi\in\Phi}|\widehat{\varphi}|^2((A^*)^{-1}\xi)-\sum_{\varphi\in\Phi}|\widehat{\varphi}|^2(\xi)=\sum_{\psi\in\Psi}|\widehat{\psi}|^2(\xi),
\end{equation}
Then take $\lambda=1$ and $\lambda=i$ and substract the equalities: 
\begin{equation}\label{eq1_1_4}
-\sum_{\varphi\in\Phi}\widehat{\varphi}(\xi)\overline{\widehat{\varphi}}(\xi+2s\pi)=\sum_{\psi\in\Psi}\widehat{\psi}(\xi)\overline{\widehat{\psi}}(\xi+2s\pi)
\end{equation}
for all $s\in\Zn\setminus A^*\Zn$. 
\par
Now take $s\in A^*\Zn$. Then
$$D_d^*(\delta_0+\lambda\delta_s)=\left\{
\begin{array}{ccc}
0&\mbox{if}&d\neq0\\
\delta_0+\lambda\delta_{(A^*)^{-1}s}&\mbox{if}&d=0.
\end{array}
\right.$$
$$\tau_{V_1,\delta_0+\lambda\delta_s}(\xi)=\tau_{V_0,\delta_0+\lambda\delta_{(A^*)^{-1}s}}((A^*)^{-1}\xi)$$
$$=\sum_{\varphi\in\Phi}|\widehat{\varphi}((A^*)^{-1}\xi)+\overline{\lambda}\widehat{\varphi}((A^*)^{-1}\xi+2\pi(A^*)^{-1}s)|^2
$$
Therefore,
$$\tau_{W_0,\delta_0+\lambda\delta_s}(\xi)=\sum_{\varphi\in\Phi}|\widehat{\varphi}((A^*)^{-1}(\xi)+\overline{\lambda}\widehat{\varphi}((A^*)^{-1}(\xi+2s\pi))|^2
-\sum_{\varphi\in\Phi}|\widehat{\varphi}(\xi)+\overline{\lambda}\widehat{\varphi}(\xi+2s\pi)|^2.$$
With (\ref{eq1_1_3}) it follows that 
$$\sum_{\varphi\in\Phi}\widehat{\varphi}((A^*)^{-1}\xi)\overline{\widehat{\varphi}}((A^*)^{-1}(\xi+2s\pi))-
\sum_{\varphi\in\Phi}\widehat{\varphi}(\xi)\overline{\widehat{\varphi}}(\xi+2s\pi)
=\sum_{\psi\in\Psi}\widehat{\psi}(\xi)\overline{\widehat{\psi}}(\xi+2s\pi).$$
The converse follows by retracing the calculations and using theorem \ref{th00_5_1}.
\end{proof}

\par
For the next theorem we assume the multiwavelet is given for a fixed GMRA and show that if some functions satisfy the equations 
discovered in theorem \ref{th1_1}, then they will be multiscaling functions associated to this GMRA.
\begin{theorem}\label{th1_2}
Let $(V_n)_{n\in\mathbb{Z}}$ be a GMRA with $\operatorname*{dim}_{V_0}(\xi)<\infty$ for a set of positive measure and 
let $\Psi$ be a NTF generator for $W_0$. The following affirmations are equivalent
\begin{enumerate}
\item
$\Phi$ is a NTF generator for $V_0$. 
\item
The equations (\ref{eq1_1_1}) and (\ref{eq1_1_2}) hold and
\begin{equation}\label{eq1_2_1}
\lim_{j\rightarrow\infty}\sum_{\varphi\in\Phi}|\widehat{\varphi} |^2((A^*)^j\xi)=0,\mbox{ for a.e. }\xi\in\Rn.
\end{equation}
\end{enumerate}
\end{theorem}

\begin{proof}
We know that (i) implies (\ref{eq1_1_1}) and (\ref{eq1_1_2}). Let's check (\ref{eq1_2_1}). By theorem 4.1 in 
\cite{BoRz}
$$\sum_{j=-\infty}^\infty\sigma_{W_0}((A^*)^j\xi)=1,\mbox{ for a.e }\xi\in\Rn$$
and
$$\sigma_{V_0}(\xi)=\sum_{j\geq 1}^\infty\sigma_{W_0}((A^*)^j\xi),\mbox{ for a.e. }\xi\in\Rn.$$
Then
$$\sigma_{V_0}((A^*)^J\xi)=\sum_{j\geq J+1}^{\infty}\sigma_{W_0}((A^*)^j\xi)\rightarrow 0,\mbox{ as }J\rightarrow\infty\mbox{ for a.e. }\xi\in\Rn.$$
But
$$\sigma_{V_0}(\xi)=\tau_{V_0,\delta_0}(\xi)=\sum_{\varphi\in\Phi}|\widehat{\varphi}|^2(\xi),$$
and (\ref{eq1_2_1}) follows immediately.
\par
For the converse we use theorem \ref{th00_9}. For all $j\geq1$, using (\ref{eq1_1_2}) for $\xi=(A^*)^j\xi$ and $s=(A^*)^js$,
$$\sum_{\psi\in\Psi}\widehat{\psi}((A^*)^j\xi)\overline{\widehat{\psi}}((A^*)^j(\xi+2s\pi))$$
$$=\sum_{\varphi\in\Phi}\widehat{\varphi}((A^*)^{j-1}\xi)\overline{\widehat{\varphi}}((A^*)^{j-1}(\xi+2s\pi))
-\sum_{\varphi\in\Phi}\widehat{\varphi}((A^*)^j\xi)\overline{\widehat{\varphi}}((A^*)^j(\xi+2s\pi)).$$
Now sum for $j\in\{1,...,J\}$. 
$$\sum_{j=1}^J\sum_{\psi\in\Psi}\widehat{\psi}((A^*)^j\xi)\overline{\widehat{\psi}}((A^*)^j(\xi+2s\pi))$$
$$=\sum_{\varphi\in\Phi}\widehat{\varphi}(\xi)\overline{\widehat{\varphi}}(\xi+2s\pi)
-\sum_{\varphi\in\Phi}\widehat{\varphi}((A^*)^J\xi)\overline{\widehat{\varphi}}((A^*)^J(\xi+2s\pi)).$$
But, by (\ref{eq1_2_1}),
$$|\sum_{\varphi\in\Phi}\widehat{\varphi}((A^*)^J\xi)\overline{\widehat{\varphi}}((A^*)^J(\xi+2s\pi))|$$
$$\leq\left(\sum_{\varphi\in\Phi}|\widehat{\varphi}|^2((A^*)^J\xi)\right)^{1/2}
\left(\sum_{\varphi\in\Phi}|\widehat{\varphi}|^2((A^*)^J(\xi+2s\pi))\right)^{1/2}\rightarrow 0,\mbox{ for a.e. }\xi.$$
as $J\rightarrow\infty$.
Hence
$$\sum_{j=1}^\infty\sum_{\psi\in\Psi}\widehat{\psi}((A^*)^j\xi)\overline{\widehat{\psi}}((A^*)^j(\xi+2s\pi))
=\sum_{\varphi\in\Phi}\widehat{\varphi}(\xi)\overline{\widehat{\varphi}}(\xi+2s\pi)$$
and, with theorem \ref{th00_9}, we can conclude that $\Phi$ is a NTF generator for $V_0$. 
\end{proof}

\begin{proposition}\label{prop1_3}
Let $V_0$ be a refinable space i.e. $V_0\subset D_AV_0$ and let $\Phi$ be a NTF generator for $V_0$. Denote by 
$V_j=D_A^jV_0$, $j\in\mathbb{Z}$. The following affirmations are equivalent:
\begin{enumerate}
\item
$$\overline{\bigcup_{j\in\mathbb{Z}}V_j}=\ltworn;$$
\item
$$\lim_{j\rightarrow\infty}\sum_{\varphi\in\Phi}|\widehat{\varphi}|^2((A^*)^{-j}\xi)=1,\mbox{ for a.e. }\xi\in\Rn.$$
\end{enumerate}
\end{proposition}

\begin{proof}
Let 
$$V:=\overline{\bigcup_{j\in\mathbb{Z}}V_j}.$$
Then, for a.e $\xi$,
$$\tau_{V,\delta_0}(\xi)=\lim_{j\rightarrow\infty}\tau_{V_j,\delta_0}(\xi).$$
But, according to proposition \ref{prop00_7},
$$\tau_{V_j,\delta_0}(\xi)=\tau_{V,\delta_0}((A^*)^{-j}\xi)=\sum_{\varphi\in\Phi}|\widehat{\varphi}|^2((A^*)^{-j}\xi).$$
If (i) holds then $\tau_{V,\delta_0}(\xi)=1$ so (ii) is immediate. 
\par
If (ii) holds then $\tau_{V,\delta_0}=1$ which implies $\delta_0\in J_{per}(\xi)$ for a.e. $\xi$, $J_{per}$ being the 
periodic range function associated to $V$. By periodicity,
$\delta_k=\lambda(-k)^*\delta_0\in\lambda(-k)^*J_{per}(\xi+2k\pi)=J_{per}(\xi)$ so that 
$\delta_k\in J_{per}(\xi)$ for all $k$ for a.e. $\xi$. This means that $J_{per}(\xi)=\ltwozn$ almost everywhere 
so $V=\ltworn$. 
\end{proof}

\begin{proposition}\label{prop1_4}
Let $V_0$ be a refinable space then 
$$(\sum_{\varphi\in\Phi}|\widehat{\varphi}|^2((A^*)^{-j}\xi))_{j\geq0}$$ is an increasing sequence for a.e. $\xi$. 
\end{proposition}

\begin{proof}
Indeed $\tau_{V_j,\delta_0}(\xi)=\sum_{\varphi\in\Phi}|\widehat{\varphi}|^2((A^*)^{-j}\xi)$ and the rest follows by the 
monotony of the local trace function.
\end{proof}

\begin{theorem}\label{th1_5}
Let $\Phi,\Psi$ be two subsets of $\ltworn$ with $\Psi$ finite and $\Phi$ countable. Suppose the following relations are 
satisfied:
\begin{equation}\label{eq1_5_0}
\sum_{\varphi\in\Phi}|\widehat{\varphi}|^2(\xi)<\infty,\mbox{ for a.e. }\xi
\end{equation}
\begin{equation}\label{eq1_5_1}
-\sum_{\varphi\in\Phi}\widehat{\varphi}(\xi)\overline{\widehat{\varphi}}(\xi+2s\pi)=
\sum_{\psi\in\Psi}\widehat{\psi}(\xi)\overline{\widehat{\psi}}(\xi+2s\pi),(s\in\Zn\setminus A^*\Zn);
\end{equation}

$$\sum_{\varphi\in\Phi}\widehat{\varphi}((A^*)^{-1}\xi)\overline{\widehat{\varphi}}((A^*)^{-1}(\xi+2s\pi))
-\sum_{\varphi\in\Phi}\widehat{\varphi}(\xi)\overline{\widehat{\varphi}}(\xi+2s\pi)=$$
\begin{equation}\label{eq1_5_2}
=\sum_{\psi\in\Psi}\widehat{\psi}(\xi)\overline{\widehat{\psi}}(\xi+2s\pi),\,(s\in A^*\Zn);
\end{equation}
\begin{equation}\label{eq1_5_3}
\lim_{J\rightarrow\infty}\sum_{\varphi\in\Phi}|\widehat{\varphi}|^2((A^*)^J\xi)=0,\mbox{ for a.e. }\xi;
\end{equation}
\begin{equation}\label{eq1_5_4}
\lim_{J\rightarrow\infty}\sum_{\varphi\in\Phi}|\widehat{\varphi}|^2((A^*)^{-J}\xi)=1,\mbox{ for a.e. }\xi.
\end{equation}
Then 
$$\{D_A^jT_k\psi\,|\,j\in\mathbb{Z},k\in\Zn,\psi\in\Psi\}$$
is a NTF for $\ltworn$.
\end{theorem}

\begin{proof}
We use the characterization from theorem \ref{th00_8}. For $j\geq 1$ and any $s\in\Zn$, using (\ref{eq1_5_2}) with $\xi=(A^*)^j\xi$ and 
$s=(A^*)^js$, we have
$$\sum_{\psi\in\Psi}\widehat{\psi}((A^*)^j\xi)\overline{\widehat{\psi}}((A^*)^j(\xi+2s\pi))=$$
$$=\sum_{\varphi\in\Phi}\widehat{\varphi}((A^*)^{j-1}\xi)\overline{\widehat{\varphi}}((A^*)^{j-1}(\xi+2s\pi))
-\sum_{\varphi\in\Phi}\widehat{\varphi}((A^*)^j\xi)\overline{\widehat{\varphi}}((A^*)^j(\xi+2s\pi))$$
Then, sum over $j\in\{1,...,J\}$:
$$\sum_{j=1}^J\sum_{\psi\in\Psi}\widehat{\psi}((A^*)^j\xi)\overline{\widehat{\psi}}((A^*)^j(\xi+2s\pi))=$$
$$\sum_{\varphi\in\Phi}\widehat{\varphi}(\xi)\overline{\widehat{\varphi}}(\xi+2s\pi)-
\sum_{\varphi\in\Phi}\widehat{\varphi}((A^*)^J\xi)\overline{\widehat{\varphi}}((A^*)^J(\xi+2s\pi))$$
An application of (\ref{eq1_5_0}) and the Schwarz inequality shows that the last sum converges to 0 as $J\rightarrow\infty$, hence one can 
conclude that 
\begin{equation}\label{eq1_5_5}
\sum_{j=1}^\infty\sum_{\psi\in\Psi}\widehat{\psi}((A^*)^j\xi)\overline{\widehat{\psi}}((A^*)^j(\xi+2s\pi))=
\sum_{\varphi\in\Phi}\widehat{\varphi}(\xi)\overline{\widehat{\varphi}}(\xi+2s\pi),\,(s\in\Zn).
\end{equation}
This, added to (\ref{eq1_5_1}), yields 
$$\sum_{j=0}^\infty\sum_{\psi\in\Psi}\widehat{\psi}((A^*)^j\xi)\overline{\widehat{\psi}}((A^*)^j(\xi+2s\pi))=0,\,(s\in\Zn\setminus A^*\Zn).$$
Also, from (\ref{eq1_5_5}) with $s=0$ and $\xi=(A^*)^{-J}\xi$, we have
$$\sum_{j=-J+1}^\infty\sum_{\psi\in\Psi}|\widehat{\psi}|^2((A^*)^j\xi)=
\sum_{\varphi\in\Phi}|\widehat{\varphi}|^2((A^*)^{-J}\xi).$$
Then (\ref{eq1_5_4}) implies
$$\sum_{j=-\infty}^\infty\sum_{\psi\in\Psi}|\widehat{\psi}|^2((A^*)^j\xi)=1,\mbox{ for a.e. }\xi,$$
and the conclusion is proved with theorem \ref{th00_8}.
\end{proof}

\section{\label{construction}A construction of NTF wavelets}
We give a construction of NTF wavelets that starts from the spectral function $\sigma_{V_0}=\tau_{V_0,\delta_0}$. We construct 
two sets of functions $\Phi$ and $\Psi$ satisfying the hypotheses of theorem \ref{th1_5} and such that 
$$\sigma(\xi)=\sum_{\varphi\in\Phi}|\widehat{\varphi}|^2(\xi)=\sum_{j\geq1}\sum_{\psi\in\Psi}|\widehat{\psi}((A^*)^j\xi)|^2.$$
We will obtain a NTF multiwavelet, $\Psi$.
\par
In the sequel, we describe the construction. The starting point is a function $\sigma$ on $\Rn$ with the following properties:
\begin{equation}\label{eq2_1}
\sigma\in L^1(\mathbb{R}^n),\quad\sigma\geq0;
\end{equation}
\begin{equation}\label{eq2_2}
\sigma(A^*\xi)\leq\sigma(\xi),\mbox{ for a.e. }\xi\in\Rn;
\end{equation}
\begin{equation}\label{eq2_3}
\mbox{If }K\mbox{ is the support of }\xi\mapsto\sigma((A^*)^{-1}\xi)-\sigma(\xi)\mbox{ then }\operatorname*{Per}(\chi_{K})\mbox{ is bounded};
\end{equation}
Recall that, for $f\in L^1(\Rn)$,
$$\operatorname*{Per}(f)(\xi):=\sum_{k\in\Zn}f(\xi+2k\pi),\quad(\xi\in\Rn).$$
\begin{equation}\label{eq2_4}
\lim_{J\rightarrow\infty}\sigma((A^*)^{-J}\xi)=1,\mbox{ for a.e. }\xi\in\Rn;
\end{equation}
\begin{equation}\label{eq2_5}
\lim_{J\rightarrow\infty}\sigma((A^*)^J\xi)=0,\mbox{ for a.e. }\xi\in\Rn.
\end{equation}
In view of theorem \ref{th1_2} (ii), proposition \ref{prop1_3} and \ref{prop1_4} and since $\sigma$ should be a spectral 
function, the conditions (\ref{eq2_1}),(\ref{eq2_2}),(\ref{eq2_4}) and (\ref{eq2_5}) are natural. Condition 
(\ref{eq2_3}) will allow us to pick a finite number of $\psi$'s. 
\par
We can restate condition (\ref{eq2_3}) as follows 
\begin{equation}\label{eq2_6}
\mbox{There is a finite partition }K_1,...,K_p\mbox{ of }K\mbox{ such that }
\end{equation}
$$\mbox{ no }K_i\mbox{ contains }\xi\mbox{ and }\xi+2k\pi\mbox{ for some }k\neq 0,\xi\in\Rn.$$
This is clear because $\operatorname*{Per}(\chi_{K})(\xi)=l$ means that there are exactly $l$ points in $K$ that are congruent to 
$\xi$ modulo $2\pi$. One way to choose such a partition is by intersecting $K$ with the intervals 
$[-\pi,\pi]^n+2k\pi$. A better way is to let $p$ be the maximum of $\operatorname*{Per}(\chi_{K})$. First pick a measurable subset $K_p$ of $K$ such that $K_p$ is 
congruent modulo $2\pi$ to $\{\xi\in[-\pi,\pi]^n\,|\,\operatorname*{Per}(\chi_{K})(\xi)=p\}$, then pick a measurable subset $K_{p-1}$ of 
$K\setminus K_p$ such that $K_{p-1}$ is congruent to 
$\{\xi\in[-\pi,\pi]^n\,|\,\operatorname*{Per}(\chi_{K})(\xi)\geq p-1\}$ and so on, finally $K_1$ is a subset of $K\setminus\cup_{i=2}^pK_i$ which is 
congruent modulo $2\pi$ to 
$\{\xi\in[-\pi,\pi]^n\,|\,\operatorname*{Per}(\chi_{K})(\xi)\geq 1\}$.
\par
Consider that we have built this partition. Define the measurable functions $\varphi_k$ for $k\in\Zn$ as follows
$$\widehat{\varphi}_k(\xi)=\left\{
\begin{array}{ccc}
\sqrt{\sigma(\xi)}&\mbox{if}&\xi\in [-\pi,\pi]^n+2k\pi\\
0&,&\mbox{otherwise}.
\end{array}
\right.
$$
Clearly then $\varphi_k\in\ltworn$ and
\begin{equation}\label{eq2_7}
\sum_{k\in\Zn}|\widehat{\varphi}_k|^2(\xi)=\sigma(\xi),\mbox{ for a.e. }\xi,
\end{equation}
and
\begin{equation}\label{eq2_8}
\widehat{\varphi}_k(\xi)\overline{\widehat{\varphi}}_k(\xi+2s\pi)=0,\quad(\xi\in\Rn,s\in\Zn\setminus\{0\},k\in\Zn).
\end{equation}
Now we construct the wavelets: for $i\in\{1,...,p\}$ define a measurable $\psi_i$ such that 
$$|\widehat{\psi}_i|^2(\xi)=\left\{
\begin{array}{ccc}
\sigma((A^*)^{-1}(\xi))-\sigma(\xi)&\mbox{if}&\xi\in K_i\\
0&,&\mbox{otherwise}.
\end{array}
\right.
$$
This is possible because (\ref{eq2_2}) holds. Discard those $\psi_i$'s which are identically 0. 
Then, of course, $\psi_i\in\ltworn$,
\begin{equation}\label{eq2_9}
\sum_{i=1}^p|\widehat{\psi}_i|^2(\xi)=\sigma((A^*)^{-1}\xi)-\sigma(\xi),\mbox{ for a.e. }\xi\in\Rn,
\end{equation}
and
\begin{equation}\label{eq2_10}
\widehat{\psi}_i(\xi)\overline{\widehat{\psi}}(\xi+2s\pi)=0,\quad(\xi\in\Rn,s\in\Zn\setminus\{0\},i\in\{1,...,p\}).
\end{equation}
(\ref{eq2_7}),(\ref{eq2_8}),(\ref{eq2_9}),(\ref{eq2_10}) and theorem \ref{th1_5} imply the fact that $\{\psi_i\,|\, i\in\{1,...,p\}\}$ is 
a NTF wavelet.

\begin{example}\label{ex2_1}
{\bf [NTF multi-wavelet sets]}
Any NTF multi-wavelet set can be obtained with this construction. Recall that a NTF multi-wavelet set is a 
NTF multiwavelet $\Psi=\{\psi_1,...,\psi_p\}$ such that each $\widehat{\psi}_i$ is the characteristic function of some measurable 
set $E_i$. 
\par
First, we give a theorem which characterizes NTF multi-wavelet sets. For a different proof when $p=1$ and some related topics see
\cite{DDGH}.
\begin{theorem}\label{th2_2}
If $\Psi=\{\psi_1,...,\psi_p\}$ and $\widehat{\psi}_i=\chi_{E_i}, (i\in\{1,...,p\})$, then
$\Psi$ is a multiwavelet set if and only if the following conditions are satisfied:
\begin{equation}\label{eq2_2_1}
E_1,...,E_p\mbox{ are mutually disjoint};
\end{equation}
\begin{equation}\label{eq2_2_2}
E_i\cap(E_i+2k\pi)=\emptyset,\quad(k\neq 0,i\in\{1,...,p\});
\end{equation}
\begin{equation}\label{eq2_2_3}
\{(A^*)^j(\cup_{i=1}^pE_i)\,|\, j\in\mathbb{Z}\}\mbox{ is a partition of }\Rn.
\end{equation}
\end{theorem}
\begin{proof}
By theorem \ref{th00_8}, $\Psi$ is a NTF multi-wavelet iff the equations (\ref{eq00_8_1}) and (\ref{eq00_8_2}) hold. Since
$\widehat{\psi}_i=\chi_{E_i}$, (\ref{eq00_8_2}) is equivalent to 
$$\chi_{E_i}((A^*)^j\xi)\chi_{E_i}((A^*)^j(\xi+2s\pi))=0,\quad(\xi\in\Rn,i\in\{1,...,p\},j\geq0,s\in\Zn\setminus A^*\Zn),$$
and, as any number can be written as $(A^*)^j\xi$ and any $k\neq0$ can be written as $(A^*)^js$, this is true iff 
(\ref{eq2_2_2}) holds. 
\par
(\ref{eq00_8_1}) rewrites as 
$$\sum_{i=1}^p\sum_{j\in\mathbb{Z}}\chi_{(A^*)^jE_i}=1,$$
which is equivalent to (\ref{eq2_2_1}) and (\ref{eq2_2_3}).
\end{proof}
Now let's see how any NTF multi-wavelet set is obtained with our construction. 
\par
Let $\Psi:=\{\psi_1,...,\psi_p\}$, $\widehat{\psi}_i=\chi_{E_i},(i\in\{1,...,p\})$ be the NTF multi-wavelet set. Define
$$\sigma(\xi)=\sum_{i=1}^p\sum_{j\geq1}\chi_{E_i}((A^*)^j\xi),\quad(\xi\in\Rn).$$
Then, theorem \ref{th2_2}, implies that $\sigma=\chi_E$, where
$$E=\bigcup_{i=1}^p\bigcup_{j\geq1}(A^*)^{-j}E_i.$$
Therefore,
\begin{equation}\label{eq2_1_1}
\sigma((A^*)^{-1}\xi)-\sigma(\xi)=\chi_{\cup_{i=1}^pE_i}(\xi),\quad(\xi\in\Rn).
\end{equation}
A simple check shows that $\sigma$ satisfies the equations (\ref{eq2_1})-(\ref{eq2_5}). 
\par
Then, we can choose the partition $K_i$ to be $K_i=E_i$ and $\widehat{\psi}_i=E_i$, $(i\in\{1,...,p\}$ so that, after the 
construction we get back our NTF multi-wavelet set $\Psi$.
\par
But how must $\sigma$ be chosen, if we want to obtain a multi-wavelet set after the construction? We saw that $\sigma$ must be 
a characteristic function of some measurable set 
$$\sigma=\chi_E.$$
The conditions (\ref{eq2_1})-(\ref{eq2_5}) can be reformulated as:
\begin{equation}\label{eq2_1_2}
E\mbox{ has finite measure };
\end{equation}
\begin{equation}\label{eq2_1_3}
E\subset A^*E;
\end{equation}
\begin{equation}\label{eq2_1_4}
\operatorname*{Per}(\chi_{A^*E\setminus E})\mbox{ is bounded};
\end{equation}
\begin{equation}\label{eq2_1_5}
\mbox{For a.e. }\xi\mbox{ there exists a }J_0\mbox{ such that }(A^*)^{-j}\xi\in E\mbox{ for }j\geq J_0;
\end{equation}
\begin{equation}\label{eq2_1_6}
\mbox{For a.e. }\xi\mbox{ there exists a }J_0\mbox{ such that }(A^*)^j\xi\not\in E\mbox{ for }j\geq J_0.
\end{equation}
The proceed with the construction and the $\widehat{\psi}_i$'s can be chosen characteristic functions. 
\par
Hence starting with $\sigma=\chi_E$ that verifies (\ref{eq2_1_2})-(\ref{eq2_1_6}), the construction yields
multi-wavelet sets. 
\par
A simple example of such a function $\sigma$, for $\ltwor$ and $A=2$, would be $\sigma:=\chi_{(a,b)}$, where the 
interval $(a,b)$ contains $0$. Then 
$$\sigma(2^{-1}\xi)-\sigma(\xi)=\chi_{(2a,a]\cup[b,2b)}(\xi),$$
and the NTF multi-wavelet set is obtained taking 
$$\widehat{\psi}_i:=\chi_{((2a,a]\cup[b,2b))\cap K_i},$$
where $K_i$ is a partition of $(2a,a]\cup[b,2b)$ as described in the construction. 
\par
To obtain single NTF wavelet sets (i.e. $p=1$), one has to start with $\sigma=\chi_E$, where $E$ verifies (\ref{eq2_1_2}),
(\ref{eq2_1_3}),(\ref{eq2_1_5}),(\ref{eq2_1_6}), and (\ref{eq2_1_4})  is replaced by 
\begin{equation}\label{eq2_1_7}
\operatorname*{Per}(\chi_{A^*E\setminus E})\leq 1.
\end{equation}
\par
If we analyse the argument before, we see that any single NTF wavelet set comes from such a construction. 
\par
For single orthonormal wavelet sets, we have the same conditions, only (\ref{eq2_1_7}) must be replaced by 
\begin{equation}\label{eq2_1_8}
\operatorname*{Per}(\chi_{A^*E\setminus E})=1.
\end{equation}
\end{example}

\begin{example}\label{ex2_2}Let $a,b>0$. Define the piecewise linear function
$$\sigma(\xi)=\left\{
\begin{array}{ccc}
\frac{1}{a}\xi+1,&\mbox{if}&\xi\in(-a,0]\\
-\frac{1}{b}\xi+1,&\mbox{if}&\xi\in[0,b)\\
0,& &\mbox{otherwise.}
\end{array}\right.$$
Then,
$$\sigma(2^{-1}\xi)=\left\{
\begin{array}{ccc}
\frac{1}{2a}\xi+1,&\mbox{if}&\xi\in(-2a,0]\\
-\frac{1}{2b}\xi+1,&\mbox{if}&\xi\in[0,2b)\\
0,& &\mbox{otherwise.}
\end{array}\right.$$
and 
$$\sigma(2^{-1}\xi)-\sigma(\xi)=\left\{
\begin{array}{ccc}
\frac{1}{2a}\xi+1,&\mbox{if}&\xi\in(-2a,a]\\
-\frac{1}{2a}\xi,&\mbox{if}&\xi\in(a,0]\\
\frac{1}{2b}\xi,&\mbox{if}&\xi\in[0,b)\\
-\frac{1}{2b}\xi+1,&\mbox{if}&\xi\in[b,2b)\\
0,& &\mbox{otherwise.}
\end{array}\right.$$
A simple check shows that $\sigma$ satisfies the conditions (\ref{eq2_1})-(\ref{eq2_5}).
Then  we can construct a NTF multiwavelet as before, we only need the partition $K_i$. For example, take 
$K_l=[(2l-1)\pi,(2l+1)\pi)$ and intersect with $(-2a,2b)$. Then, let 
$$\eta(\xi):=
\left\{
\begin{array}{ccc}
\sqrt{\frac{1}{2a}\xi+1},&\mbox{if}&\xi\in(-2a,a]\\
\sqrt{-\frac{1}{2a}\xi},&\mbox{if}&\xi\in(a,0]\\
\sqrt{\frac{1}{2b}\xi},&\mbox{if}&\xi\in[0,b)\\
\sqrt{-\frac{1}{2b}\xi+1},&\mbox{if}&\xi\in[b,2b)\\
0,& &\mbox{otherwise.}
\end{array}\right.$$
and $\widehat{\psi}_l=\eta\chi_{K_l}$ for $l\in\mathbb{Z}$ with the property that $[(2l-1)\pi,(2l+1)\pi)\cap(-2a,2b)\neq\emptyset$.
Then $\{\psi_l\}_l$ is a NTF multiwavelet. 
\par
If $a+b\leq\pi$ then we need only one $K_i$ and we can let $K_1=(-2a,2b)$ so that in this case $\eta$ is a NTF wavelet. 
\end{example}

\begin{remark}\label{rem2_3}
The NTF multi-wavelet sets are always semi-orthogonal, that is, if 
$$W_j=\overline{\mbox{span}}\{D_A^jT_k\psi\,|\, k\in\Zn,\psi\in\Psi\},\quad(j\in\mathbb{Z}),$$
then $W_j\perp W_{j'}$ for $j\neq j'$. \par
This can be seen from the fact that $\widehat{D}_A^j\widehat{T}_k\widehat{\psi}$ and 
$\widehat{D}_A^{j'}\widehat{T}_{k'}\widehat{\psi}'$ are disjointly supported for $j\neq j'$ and any
$k,k'\in\Zn$, $\psi,\psi'\in\Psi$ (this is guaranteed by (\ref{eq2_2_1}) and (\ref{eq2_2_3})).
\par
On the other hand, the NTF multiwavelet in example \ref{ex2_1} is not semi-orthohgonal because $\eta$ and $\widehat{D}_A\eta$ are 
positive functions that have an overlap in their supports.

\end{remark}

\end{document}